\documentclass[11pt]{article}
\usepackage[T1]{fontenc}
\usepackage[french,english]{babel}

\usepackage{graphicx}
\DeclareGraphicsExtensions{.jpg}
\makeatletter
\def\@begintheorem#1#2{\trivlist
      \item[\hskip \labelsep{\bf #1\ #2}]\rm}
\def\@opargbegintheorem#1#2#3{\trivlist
      \item[\hskip \labelsep{\bf #1\ #2\ (#3)}]\rm}
\makeatother         
\def\rit{\hbox{\it I\hskip -2pt R}}

\title{Robert de Montessus de Ballore's 1902 theorem on algebraic continued fractions :
 genesis and circulation}
\author{Herv\'e Le Ferrand \footnote{Universit\'e de Bourgogne, Institut de Math\'ematiques de Bourgogne, 9 avenue Alain Savary, BP 47870, 21078 Dijon, France, email: leferran@u-bourgogne.fr}}
\date{\today}

\begin{document}
\maketitle
\markright{Circulation of a theorem}


{\bf Abstract}
Robert de Montessus de Ballore proved in 1902 his famous theorem on the convergence of Pad\'e approximants of meromorphic functions. In this paper, we will first describe the genesis of the theorem, then investigate its circulation. A number of letters addressed to Robert de Montessus by different mathematicians will be quoted to help determining the scientific context and the steps that led to the result. 
In particular, excerpts of the correspondence with Henri Pad\'e in the years 1901-1902 played a leading role. 
The large number of authors who mentioned the theorem soon after its derivation, for instance N\"orlund and Perron among others, indicates a fast circulation due to factors that will be thoroughly explained.

{\bf Key words}
Robert de Montessus, circulation of a theorem, algebraic continued fractions, Pad\'e's approximants.

{\bf MSC} : 01A55 ; 01A60


\section{Introduction}

\paragraph{}
This paper aims to study the genesis and circulation of the theorem on convergence of {\it algebraic continued fractions} proven by the French mathematician Robert de Montessus de Ballore (1870-1937) in 1902. The main issue is the following : which factors played a role in the elaboration then the use of this new result ? Inspired by the study of Sturm's theorem by Hourya Sinaceur \cite{Sina}, the scientific context of Robert de Montessus' research will be described. Additionally, the correlation with the other topics on which he worked will be highlighted, as well as the major points that led him to the result. In \cite{Gold_1}, Catherine  Goldstein raised the following issues for any research in history of mathematics :
  \og{\it What is the historical description of a theorem in mathematics ? } \fg  ; \og{\it How was the theorem read and understood ? How was it used ?} \fg and : \og{\it understanding a paper, historical as well as mathematical, means determining on which knowledge it relies ;  on which options, for what public;  which context is required to interpret the text ; which points are needed to build an answer.} \fg\  In the case of Robert de Montessus, the challenge is to find the mathematical approaches of the theorem and the reasons for other mathematicians to refer to his work.

Robert de Montessus \cite{Lef}, aged 32, began his thesis in 1902 at the {\it Facult\'e des Sciences of Paris}, under the supervision of Paul Appell. His defense\footnote{His research on algebraic continued fractions was awarded by the Grand Prix de l'Acad\'emie des Sciences Math\'ematiques in 1906. The price was shared with Robert de Montessus, Henri Pad\'e and Andr\'e Auric. } on the 8th of Mai 1905\footnote{Paul Appell (chairman), Henri Poincar\'e and Edouard Goursat were present in the jury.}  led to a publication in the Bulletin de la Soci\'et\'e Math\'ematique de France (pages 28 to 36, volume 10), where he raised these two issues :

\begin{quotation}
\begin{it}

Does a sequence of consecutive convergents from a table constituted by normal fractions define a function similar to the function defined by the series that has led to the table ? And, if yes, does the continuous fraction corresponding to the progression extend the series outside its circle of convergence ?

\end{it}
\end{quotation}

Robert de Montessus gave a positive answer to these questions in the particular case of the following function : a meromorphic function, analytic around zero, and for which one consider as the sequence of rational fractions, the elements in a row of the Pad\'e table\footnote{A definition of this table will be given in the section 2.1.} associated to the function .


During the last forty years, Robert de Montessus de Ballore has been quoted a large number of times in the mathematic field in connection with the approximation using rational fractions\footnote{With one or several variables.} or with continued fractions\footnote{These two topics are part of what is called Rational Approximation.}. His name appears in the statement of a new result or as a reference\footnote{Inside the article or in the bibliography.}  to the theorem of 1902. As an example, E.B. Saff published in 1972 {\it An extansion of Montessus de Ballore's theorem on the convergence of interpoling rational functions} \cite{Saff_1}. This article\footnote{E.B Saff dedicated this article to his thesis supervisor J.L. Walsh. The latter quoted Robert de Montessus in his book published in 1935, entitled {\it Interpolation and Approximation by Rational Functions in the Complex Domain}.}  states the theorem using Pad\'e approximants and the bibliography includes a reference to Montessus' article published in 1902 in the Bulletin de la Soci\'et\'e Math\'ematique de France. Moreover, Annie Cuyt and Doron Lubinsky published in 1997 {\it A de Montessus theorem for multivariate homogeneous Pad\'e approximants} \cite{Cuyt}. Here, the name of {\it Montessus de Ballore} appears several times through a less specific book : it is the work of Perron \cite{Perron} on continued fractions that is mentioned as a reference, even if Montessus' theorem is cited. More recently in 2011, on can find the theorem stated with a direct reference to the article of 1902 in {\it From $QD$ to $LR$ and $QR$, or, how were the $QD$ and $LR$ algorithms discovered?} \cite{Gut}, by M. H. Gutknecht M.H. and B.N. Parlett.

In which context was the theorem of 1902 elaborated ? Who else was doing research on the convergence of algebraic continued fractions and what was the goal at that time ? These are the challenging issues we will try to answer in the first part. To make it easier for the reader, we will explain at the beginning of this section some definitions and results on continued fractions and Pad\'e approximants\footnote{In 1902, the term of Pad\'e approximants was not used yet. Instead, it was {\it fractions approch\'ees}.}. In the second part, the circulation of the theorem will be studied : which authors cited Robert de Montessus ? When ? What did Robert de Montessus to further improve the dissemination of his work ?


Methodologically, we proceeded as following : we used three types of sources. The books of C. Brezinski \cite{Brez1} and \cite{Brez2}\footnote{This book on the history of continued fractions by C. Brezinski is an inexhaustible source of information. Indeed, it gives very thorough mathematical explanations and the bibliography included is quite complete.} were the starting point. Then we also find a lot of valuable informations in electronic sources and in letters\footnote{Letters and documents that belonged to Robert de Montessus have been very kindly found thanks to his family. They are now in the archives of the Universit\'e Pierre et Marie Curie in Paris, France.}  received by Robert de Montessus from other mathematicians working on topics close to algebraic continued fractions. These letters represent a crucial and a new source of information on the process of elaboration of the new theorem and its circulation. To illustrate this, excerpts will be given. Moreover, for the study of the circulation\footnote{Circulation or dissemination are used here with the same meaning.}  of the theorem, we linked different electronic databases, such as Jahrbuch \cite{site2}, Numdam \cite{site1}, Jstor \cite{site3}, Gallica \cite{site4}, Internet Archive \cite{site6} and G\"ottinger Digitalisierungszentrum \cite{site7}. The goal is to gather, with the most exhaustivity possible, all the mathematical articles that mention in any way the theorem of Robert de Montessus. The book written by C. Brezinski \cite{Brez1} on Henri Pad\'e's work, will enable us to position Robert de Montessus' work as regard to Henri Pad\'e's.

\section{Robert de Montessus and continued fractions}

\subsection{Mathematical objects: continued fractions, algebraic continued fractions, Pad\'e approximants}
The mathematical objects used by Robert de Montessus will be described the same way they were at the beginning of the 20\up{th} century.
First, a continued fraction is given by two sequences of natural integers, real or complex numbers $q_{0},q_{1},q_{2},\ldots$, $p_{1},p_{2}$, and by a calculation procedure   :


$$
q_{0}+\frac{p_{1}}{q_{1}+\frac{p_{2}}{q_{2}+\frac{p_{3}}{q_{3}+\cdots}}}
$$
For example, starting with a fraction and with the use of Euclide's algorithm, one gets :
$$
\frac{105}{24}=4+\frac{1}{2+\frac{1}{1+\frac{1}{2}}}.
$$

A similar procedure but using the floor function gives for the number\footnote{In {\it Additions au m\'emoire sur la r\'esolution des \'equations num\'eriques}, published in 1770, Lagrange proved that every {\it quadratic} number, like $\sqrt{3}$, can be rewritten as a continuous and periodic fraction.}   $\sqrt{3}$:

$$
\sqrt{3}=1+(\sqrt{3}-1)=1+\frac{1}{\frac{1}{\sqrt{3}-1}}=1+\frac{1}{\frac{\sqrt{3}+1}{2}}=
1+\frac{1}{1+\frac{1}{2+\frac{1}{1+\frac{1}{2+\cdots}}}}
$$

The convergent of order $k$ of the continued fraction is the usual fraction :
$$
\frac{A_{k}}{B_{k}}=q_{0}+\frac{p_{1}}{q_{1}+\frac{p_{2}}{q_{2}+\cdots+\frac{p_{k}}{q_{k}}}}.
$$

The convergents give in certain circumstances a good approximations of the value of the continued fraction. As an illustration, in the17\up{th} century, the Dutch astronomer Christian Huygens used continued fractions and their convergents to build planetary automatons\footnote{To fix the number of gear teeth of the notched wheels of his machines, Huygens faced fractions with huge numbers, namely astronomical distances. He could approximate in a satisfactory way these quantities by quotients with small integer, using convergents of the developments in continued fractions.} in {\it Descriptio automati planetarii} \cite{Huygens}. 

The following recurrence relations link the numerators and the denominators of the convergents :
$$
A_{k}=q_{k}A_{k-1}+p_{k}A_{k-2}\ ;\ B_{k}=q_{k}B_{k-1}+p_{k}B_{k-2}.
$$
with first terms:
$$
A_{-1}=1,\ A_{0}=q_{0},\ A_{1}=q_{1}q_{0}+p_{1}
$$
and
$$
B_{-1},\ B_{0}=1,\ B_{1}=q_{1}.
$$

These two equations show in the first place a relation between the continued fractions and the {\it equations with finite linear differences}. If in \cite{Brez2} there are different references on the relation between continued fractions and equations with differences, we will mention here the thesis of Niels N\"orlund entitled {\it Fractions continues et diff\'erences r\'eciproques} and published in 1911 in  Acta Mathematica. Indeed, as we will see later on, N\"orlund contacted Robert de Montessus in 1910 about his work on algebraic continued functions. In the second place and with certain hypotheses, it is possible to write the $p_{k}$ and $q_{k}$ with the quantities $A_{i}$ and $B_{j}$. It is thus possible to build continued fractions whose convergents are a sequence of given fractions.

Furthermore, if the sequences  $p_{i}$ and $q_{i}$ are sequences of polynomials, an {\it algebraic continued fraction} has a similar definition. For example, Heinrich Lambert obtained in 1768 a development of $\tan(x)$ \footnote{The result can be found in {\it M\'emoire sur quelques propri\'et\'es remarquables des quantit\'es transcendantes, circulaires et logarithmiques}, Acad\'emie des Sciences of Berlin, 1768. Besides, H. Lambert used the development of  $\tan$ to prove $\pi$ was irrational.}, without proving the convergence\footnote{There has been a debate concerning this. It will be mentioned later when we will talk about Alfred Pringsheim.} :

$$
\tan(x)=\frac{x}{1-\frac{x^{2}}{3-\frac{x^{2}}{5-\frac{x^{2}}{7-\frac{x^{2}}{9-\cdots}}}}}
$$  
This development in algebraic continued fraction can be obtained with a software performing algebraic computation such as Maple. Thus, using {\it cfrac} in Maple, one gets :
$$
exp(z)=1+\frac{z}{1+\frac{z}{-2+\frac{z}{-3+\frac{z}{2+\cdots}}}}
$$
This equation is an example of a development in continued fraction that is part of the theory of Pad\'e approximants\footnote{Henri Pad\'e used in his thesis \cite{Pad\'e_1} the terms of {\it fraction continue simple}, pages 42-93}. These convergents namely correspond to a certain displacement in Pad\'e's table associated with the exponential function. But what is the a  Pad\'e approximant of a function ? \\
Consider $f(z)=\sum_{i=0}^{+\infty}a_{i}z^{i}$ a function\footnote {We use $z$ for the variable. Instead, Henri Pad\'e and Robert de Montessus used the letter $x$.} that can be developed in power series at the origin. The Pad\'e's approximant of $f$, $\left\lbrack L/M\right \rbrack(z)$, is a rational fraction with a numerator having a degree inferior to $L$ and the denominator a degree inferior to $M$, with

$$
f(z)-\left\lbrack L/M\right \rbrack(z)=O\left(z^{L+M+1}\right)
$$
which means that the development of $f$ and of $\left\lbrack L/M\right \rbrack$  are the same until the order $z^{L+M}$ is reached. For example, the approximant $\left\lbrack 3/4\right\rbrack$ of $\displaystyle{\exp(z)=1+z+\frac{z^{2}}{2!}+\frac{z^{3}}{3!}+\cdots}$ is:

\begin{eqnarray*}
\lefteqn{\frac{1+\frac{3}{7}z+\frac{1}{14}z^{2}+\frac{1}{210}z^{3}}{1-\frac{4}{7}z+\frac{1}{7}z^{2}-\frac{2}{105}z^{3}+\frac{1}{840}z^{4}} =}\\
& &   1+z+\frac{1}{2}z^{2}+\frac{1}{6}z^{3}+\frac{1}{24}z^{4}+\frac{1}{120}z^{5}+\frac{1}{720}z^{6}+\frac{1}{5040}z^{7}\\
& & +\frac{17 }{705600}z^{8}+\cdots
\end{eqnarray*}
These rational fractions are placed in a array called {\it Pad\'e table}. A good choice in this array of a sequence of aproximants\footnote{This is equivalent to choosing a displacement in the table.} generated a {\it continued fraction} whose {\it convergents} are the chosen fractions.This result can be found in the systematic study by Pad\'e on these approximants in his thesis defended in 1892. Thus, the convergents of the continued fraction of the function $exp$ that were given few lines before are its Pad\'e approximants\footnote{Besides, the software Maple gives this information.}, $\left\lbrack n/n\right\rbrack$ or $\left\lbrack n+1/n\right\rbrack$.


In the middle of the 18\up{th} century, Euler used a algebraic continued fraction to \og sum  \fg  the series 
$$ 
1!-2!+3!-\cdots
$$
(see \cite{Borel, Brez2, Vara}). 
Then during the 19\up{th} century, studies using algebraic continued fractions will be produced in a large quantity, as it is highlighted in \cite{Brez2} pages 190-259. For instance, the work of Lagrange on arithmetic and algebraic continued fractions are numerous\footnote{We advise to consult the bibliography given in \cite{Brez2} and the website Mathdoc, http://portail.mathdoc.fr/OEUVRES/, where the a link to reach the work of Lagrange that has been digitized and can be read on Gallica.}. Besides, C. Brezinski commented in \cite{Brez2} page 139 the dissertation of Lagrange called {\it Sur l'usage des fractions continues dans le calcul int\'egral}\footnote{This dissertation was read the July the 18\up{th} in 1776 by the Acad\'emie royale des Sciences et Belles Lettres of Berlin.}  and talked about {\it birth certificate Pad\'e approximants}. Indeed, Lagrange developed the function $(1+x)^{m}$ as an algebraic continued fraction. And he moreover noticed that consecutive convergents have developments in power series that are the same until a certain order\footnote{For a given convergent, this order is the sum of the degrees of the numerator and denominator.} to those of  $(1+x)^{m}$. It is Henri Pad\'e in his thesis \cite{Padethese}, who achieved a systematic study of Pad\'e approximants. He defined then studied the array formed by these fractions. He pointed out in particular the structure in blocks of the array. Then he linked the algebraic continued fractions with the Pad\'e approximants. The major question was: starting with a sequence of rational fractions, how is it possible to build an algebraic continued fraction whose convergents are the initial fractions ? Finally Pad\'e devoted his thesis to the particular case of the exponential function.  He dedicated it to Charles Hermite. We remind the reader that Charles Hermite used Pad\'e approximants in his derivation of the transcendence of $e$ in 1873\footnote{The analysis of the derivation of Hermite is very well explained by Michel Waldschmidt on Bibnum, http://bibnum.education.fr/}.

\subsection{The context}
Robert de Montessus was awarded his degree in science in 1901. He started soon after his thesis and it is thanks to Paul Appell that he began to work on the issue of the convergence of algebraic continued fractions. This has been found in letters that we will give some excerpts. This leads Robert de Montessus to get insight into the work of Henri Pad\'e, but also of Edmond Laguerre. These three mathematicians Laguerre, Appell and Pad\'e have also in common to be close to Charles Hermite. Indeed, Laguerre has graduated from the Ecole polytechnique and was a close friend of Joseph Bertrand, the brother-in-law of Hermite. Appell got Hermite as his Professor and he married one of his nieces. And Hermite supervised\footnote{There were no real thesis' supervisor at that time, so that here {\it supervising} means more {\it following}.} Pad\'e's doctoral dissertation.


Edmond Laguerre is in particular well known for the orthogonal polynomials to which he gave his name. Despite a short life -he died at 52- and health problems, he published not less than 150 articles \cite{Lag1} and \cite{Lag2}.  Henri Poincar\'e wrote in the preface of {\it Oeuvres compl\`etes de Laguerre} \cite{Lag1} :

\begin{quotation}
\begin{it}
The study of algebraic functions will without any doubt help us one day to write the functions as more convergent developments than power series ; however, few geometricians have tried to investigate in this field that will certainly give us plenty of surprises ; Laguerre went into this field thanks to his research on polynomials solving linear differential equations. He obtained many results but I want to cite only one that is the most amazing and the most suggestive. From a divergent series, a convergent continued fractions car be deduced : this is the new legitimate use of divergent series and it has a great future ahead.

\end{it}
\end{quotation}
In the same publication can be found two articles entitled {\it Sur la r\'eduction en fractions continues d'une fonction qui satisfait \`a une \'equation lin\'eaire du premier ordre \`a coefficients rationnels} that describe work that Robert de Montessus will further investigate on \cite{RMB1}. If Henri Poincar\'e underlines at the end of the 19\up{th} century the extreme importance of studying algebraic continued fractions, it is because he is working on it. First of all he published articles on continued fractions such as \cite{Poincare1, Poincare2} and then used these fractions to solve differential equations \cite{Poincare3}. Moreover, he also investigated how to sum a divergent series, like in the second volume of {\it m\'ethodes nouvelles de la m\'ecanique c\'eleste} published in 1892. Emile borel, in \cite{Borel}, dedicated a chapter to explain the theory of  {\it asymptotic series} that he attributed the creation to both Henri Poincar\'e and Thomas Stieltjes.


Paul Appell was also interested by continued fractions. One of his first article is  {\it Sur les fractions continues p\'eriodiques} \cite{Appell_1}. Nevertheless his influence on Robert de Montessus was certainly more due his researches in mathematical analysis\footnote{Appell analyzed first his work in \cite{Appell_2}. It is very likely that this note has been written for his application to the Academy of Science. Appell pointed out the importance of these development in series for the cases of the potential theory, the differential linear equations and the study of hyper geometrical functions.}, particularly his work on developments in series of functions of one or several variables


Henri Pad\'e on his side, underlined the importance of the issue of the convergence of algebraic continued fractions in {\it Recherches sur la convergence des d\'eveloppements en fractions continues d'une certaine cat\'egorie de fonctions} \cite{Pad\'e_2}. He wrote :
\begin{quotation}
\begin{it}
The question of the convergence is not mentioned in the researches on the transformation of continued fractions in power series of Lagrange and Laplace ; there are not either mentioned in the study of Gauss on the well known continued fraction resulting from the quotient of two hyper geometrical series. Nothing can be found on this topic in the huge work of Cauchy ; and, near the middle of the 19\up{th} century, if we can cite the large number of studies -not all very fruitful- of Stern and Seidel and those of Heine, we need to reach the posthumous work of Riemann (1863) : Sullo svolgimento del quoziente di due serie ipergeometrie in frazione continua infinita to find the first really relevant work. However, this work was still not complete, despite the restoration of Mr. Schwarz, and the nice and highlighting M\'emoires by Thom\'e  on the same topic and published in the volumes 66 and 67 (1866, 1867) of the Journal by Crell (...)


Since, the studies of Laguerre, Halphen, Mr. Pincherle and Markoff followed one after the other: and, for the last ten years, those of Stieltjes (1894), MM. von Koch, van Vleck and Pringsheim which are the main studies ; finally and recently, two notes to the Comptes rendus and two short dissertations published in the Bulletin de la Soci\'et\'e Math\'ematique of France (1902) and in the Annales de la Soci\'et\'e Scientifique of Bruxelles (1903), by Mr. R. de Montessus de Ballore.

\end{it}
\end{quotation}
For Henri Pad\'e, the first relevant results on the convergence of an algebraic continued function appeared at the beginning of 1860. It could be that for him, the term of {\it Oeuvres} meant a group of several works. Else his history is not complete. Why did he not mention Adrien Legendre and his proof of the convergence of the development of the function tangent in algebraic continued fraction, as it is shown above by Lambert ? This derivation can be found for example in the {\it El\'ements de G\'eom\'etrie} by Legendre \cite{Legendre} published in 1837. The works of Bernhard Riemann and of Ludwig Wilhelm Thom\'e are related to the convergence of continued fractions that are build on quotients of hypergeometrical fractions \cite{Thom}. The American mathematician Edward Burr Van Vleck defended his thesis entitled {\it Zur Kettenbruchentwickelung Hyperelliptischer und Ahnlicher Integrale} in G\"ottingen under the supervision of F\'elix Klein \cite{bioVan}, and he developed in 1903 the questions raised above as well as the studies of mathematicians mentioned by Pad\'e\footnote{ The bibliography given by Van Vleck is particularly interesting.} \cite{VanV}.  Henri Pad\'e quoted from Andrei Andreyevich Markov but it is worth pointing out both Pafnuty Lvovich Chebyshev (Tchebychev) and Konstantin Aleksandrovich Poss\'e who was as Markov a student of Chebyshev. Georges Henri Halphen studied particularly the convergence of continued fractions associated with $\sqrt{X}$, where $X$ is a polynomial of degree  $3$ with real coefficients and roots \cite{Brez2} pp 200-201. Henri Pad\'e mentioned also the Italian mathematician Salvatore Pincherle because he linked the results obtained for linear equations with finite differences with algebraic continued fractions (see {\it notice sur travaux}\cite{Pincher}). Niels Fabian Helge von Koch, student of Mittag-Leffler, generalized a result of convergence found by Thomas Stieltjes \cite{Koch}. Emile Borel devoted a chapter\footnote{Moreover in this book, Emile Borel talked about the work of Henri Pad\'e and thus about the link between Pad\'e approximants and algebraic continued fractions. } to continued fractions and to the theory of Stieltjes in his {\it Le\c cons sur les divergentes }\cite{Borel2}. at the page 63, Emile Borel wrote :

\begin{quotation}
\begin{it}
The starting point of the research of Stieltjes is the continued fraction

$$
\frac{1}{a_{1}z+\frac{1}{a_{2}+\frac{1}{a_{3}z+\frac{1}{a_{4}+\cdots}}}}
$$

where $a_{n}$ are real positive numbers and $z$ is a complex variable. 
\end{it}
\end{quotation}
then
\begin{quotation}
\begin{it}
The only case that will matter here is the one where the series $\sum a_{n}$ is divergent. The continued fraction is then convergent and defines an holomorphic function in the plane of the complex variable, points with a negative real part excepted (...)


The continued fraction can be developed in series with power of $\displaystyle{\frac{1}{z}}$ (...) This development can be then obtained as following:

$$
F=\frac{c_{0}}{z}-\frac{c_{1}}{z^{2}}+\frac{c_{2}}{z^{3}}\cdots
$$

The numbers $c_{0}, c_{1}, c_{2}, \cdots$ are positive ; Stieltjes gives a way to get their expressions in function of  $a_{n}$ but these expressions are complicated. On the contrary, the $a_{n}$ can be written using the  $c_{n}$ in very elegant form (...)

\end{it}
\end{quotation}

The work of Stieltjes on  {\it th\'eorie analytique des fractions continues}\footnote{Several authors consider Stieltjes as the father of the analytic theory of continued fractions \cite{Brez2}\cite{Jones}.} is  not separable of the moment problem,i.e. given a sequence of real numbers $(m_{n})_{n=0}^{\infty}$, does a positive measure $\mu$ over $\rit$ exists, so that for every $n$, 

$$
\displaystyle{\int_{-\infty}^{+\infty}x^{n}\mu(dx)=m_{n}} \ ? 
$$

Emile Borel further developed this point : 
\begin{quotation}
\begin{it}

Therefore, the main goal of the M\'emoire de Stieltjes is to introduce an analytical element of a defined integral that is more easy to deal with and that leads to the divergent series by its development in powers of $\frac{1}{z}$. On the other hand, for a given series, is it possible to build an integral using the continued fraction as an intermediate to connect the series and the integral (...)


Stieltjes considered the integral
$$
J=\int_{0}^{\infty}\frac{f(u)du}{z+u}
$$

where $f(u)$ is a function supposed to be nonnegative (...) 

It is easy to find its formal development. We have
$$
J=\int_{0}^{\infty}\left( \frac{1}{z}-\frac{u}{z^{2}}+\frac{u^{2}}{z^{3}}-\frac{u^{3}}{z^{4}}+\cdots\right) f(u)du
$$

So if we consider
$$
c_{n}=\int_{0}^{\infty}f(u)u^{n}du\quad (n=0, 1, 2, \ldots)
$$
we obtain
$$
F=\frac{c_{0}}{z}-\frac{c_{1}}{z^{2}}+\frac{c_{2}}{z^{3}}\cdots
$$
\end{it}
\end{quotation}

After some hypothesis and according to Emile Borel, we have:
\begin{quotation}
\begin{it}

Indeed in that case, the continued fraction deducted from the series is convergent and defines an analytical function $F(z)$ that is regular in all the plane except on the negative part of the real axis ; the equality

$$
F(z)=\int_{0}^{\infty}\frac{f(u)du}{z+u}
$$

enables us to determine the function $f(u)$.
\end{it}
\end{quotation}
Thus  Stieltjes solved the moment problem in a special case.

The German mathematician Alfred Pringsheim demonstrated several results of convergence of continued fractions \cite{Brez2}, but also studied the proofs of Lambert and Legendre on the irrationality of $\pi$ in an article published in 1901. One critic\footnote{Alfred Pringsheim explained that Lambert demonstrated the convergence of the algebraic continued fraction associated with $\tan(x)$, see the article by R Wallisser, {\it On Lambert's proof of the irrationality of $\pi$}, in Algebraic number theory and Diophantine analysis, Graz, 1998 (de Gruyter, Berlin, 2000,) pp 521-530. } can be found in the second part of the volume 25, year 1901 of the Bulletin des Sciences Math\'ematiques (pp 86-88).


Every study on the convergence of algebraic continued fractions is linked to the question of the {\it analytic continuation}. This is one of the fundamental topic in Analysis at the end of the 19\up{th} century and the beginning of the 20\up{th}. In this paper we will not give a complete history on it, but we still need to mention that Walter William Rouse Ball tackles this issue in the paragraph {\it Analyse}\footnote{This paragraph has an asterisk because it has been added by Robert de Montessus. He quoted in particular Hadamard and the Finnish mathematician Ernst Lindel\"of.} in the French version of his book of history of mathematics published in 1907 \cite{Rouse}. Robert de Montessus contributed to the writing.
Two other mathematician have corresponded with Robert de Montessus, namely Jacques Hadamard and Eug\`ene Fabry. It his well-known that Jacques Hadamard obtained important results in this field, between 1888 and 1902. His dissertation defended in 1892 is indeed entitled {\it Essai sur l'\'etude des fonctions donn\'ees par leur d\'eveloppement de Taylor} \cite{Hadam1}. We do not know how much Hadamard and Robert de Montessus have written each other \footnote{We have only found one letter, without a date on it, where Hadamard indicated to Robert de Montessus how to find his article {\it Sur certaines surfaces minima}.}. However, Robert de Montessus used the work of Jacques Hadamard for the proof of his theorem in 1902. As for the French mathematician Eug\`ene Fabry, Professor at the university of Montpellier, one letter\footnote{See the following for an excerpt of this letter.} addressed to Robert de Montessus in 1901 explained the results he got in this field. We have to cite the article of Fabry published in 1869 in the  Annales scientifiques de l'Ecole Normale Sup\'erieure,  {\it Sur les points singuliers d'une fonction donn\'ee par son d\'eveloppement en s\'erie et l'impossibilit\'e du prolongement analytique dans des cas tr\`es g\'en\'eraux}.

\subsection{Genesis}

How did Robert de Montessus start his work on algebraic continued fractions? According to the classifications in \cite{Leloup} and \cite{Gispert}, the result found by Robert de Montessus is part of the field named Analysis. However, Juliette Leloup wrote page 68 of  \cite{Leloup}:
\begin{quotation}
\begin{it}
As for the dissertation in arithmetics, these form by themselves a special corpus since the reports are always full of praise as soon as it has strong links with the theory of functions : the dissertations of Cotty, Ch\^atelet, Got and
Chapelon are such examples, the applications of the new tools of the theory of functions to the theory of continued fractions in the thesis of Montessus de Ballore is another example.

\end{it}
\end{quotation}
The author might consider that continued fractions are at the border of arithmetics and analysis.

Robert de Montessus was awarded his degree in mathematical science October the 24\up{th}, 1901\footnote{According to what is written in the statement of his dissertation}. He started afterwards a doctoral thesis under the supervision of Paul Appell. Paul Appell has been a member of the Academy of Science since 1892 and became the dean of the Faculty of Science of Paris in 1903. He supported Robert de Montessus during all his career\footnote{We are able to claim this fact according to the bunch of letters from Paul Appel found in the collection  Robert de Montessus. It is letters of recommendation, or for example, letters in which he ask for the authorization to give a lecture at the Sorbonne.}. A letter from Charles-Ange Laisant\footnote {Charles-Ange Laisant, 1841-1920, graduated from the Ecole Polytechnique in Paris, was a mathematician, French politician and founder of the {\it Interm\'ediaire des math\'ematiciens}  and with Henri Fehr of the {\it Enseignement math\'ematique}.} written on the 30/10/1900\footnote{This letter can be found in the collection of letters received by Robert de Montessus de Ballore from C. A. Laisant during the period 1897-1937, lettres de C. A. Lais}  shows that Robert de Montessus wanted to be closer to Paul Appell :

\vspace{0.8cm}
\begin{minipage}[c]{.46\linewidth}
    \begin{quotation}
{\it Cher coll\`egue,

Une fois la rentr\'ee faite, en Novembre, j'esp\`ere voir M. Appell et lui parler de vous, selon votre d\'esir. Nous vous ferons signe le moment venu [...]
}
\end{quotation}  
   \end{minipage} \hfill
   \begin{minipage}[c]{.46\linewidth}
\begin{quotation}
{\it Dear colleague,

After the starting of the classes in November, I hope to see Mr. Appell and talk to him about you as you wish. We contact you nearer the time [...]
}
\end{quotation}
   \end{minipage}



\vspace{0.8cm}
Robert de Montessus was appointed lecturer in 1902 at the catholic university of Lille, and he was then nominated  assistant professor in 1903 and gave lectures in Special mathematics and rational mechanics\footnote{ Thus Robert d'Adh\'emar, at that time assistant professor there, mentioned in the minutes of the council of the Facult\'e des Sciences de l'Universit\'e de Lille on November the 7\up{th}, 1902, that a substitute teacher was needed for the class of Special mathematics. He brought the name of Robert de Montessus. The minutes of these meetings between 1886 and 1924 are preserved int he archives of the Facult\'e des Sciences de l'Universit\'e Catholique de Lille at the number S7E.}.


The first signs of the theorem that we found appear in a letter from Henri Pad\'e sent on November the 26\up{th}, 1901\footnote{Collection of letters received by Robert de Montessus de Ballore during the period 1897-1937, from Pad\'e.} :

\vspace{0.8cm}
\begin{minipage}[c]{.46\linewidth}
\begin{quotation}
{\it

[...] J'ai d'ailleurs vu M. Appell au commencement ce mois, il m'a parl\'e de vous et m'avait fait parvenir votre lettre.

La th\'eorie des fractions continue offre un champ tr\`es vaste de recherches, mais o\`u il n'est pas toujours facile de se rendre compte \`a l' avance des difficult\'es que l'on rencontrera. C'est donc avec toutes sortes de r\'eserves que je vous indiquerai, comme devant pr\'esenter un grand int\'er\'et, une \'etude approfondie de la g\'en\'eralisation des fractions continues. Je n'ai fait qu'effleurer le sujet dans un m\'emoire qui a paru, il y a quelques ann\'ees dans le journal de M. Jordan, que je me fais un plaisir de vous envoyer en m\^eme temps que cette lettre. Vous y trouverez l'indication d'un m\'emoire de M. Hermite, sur le m\^eme sujet. Dans ce m\'emoire, M. Hermite arrive \`a la m\'ethode des polyn\^omes associ\'es par des consid\'erations [...] de calcul int\'egral : il serait, sans doute, aussi bien int\'eressant d'approfondir davantage le rapport entre le calcul int\'egral et les lois de r\'ecurrence de la th\'eorie des fractions continues.

Le r\'esultat que vous m'annoncez avoir obtenu me para\^it des plus remarquables, mais doit \^etre soumis \`a des exceptions assez nombreuses. Je lirai avec plaisir votre d\'emonstration quand vous l'aurez publi\'ee.
[...] }
\end{quotation}
   \end{minipage}
\begin{minipage}[c]{.46\linewidth}
\begin{quotation}
{\it [...] Moreover, I have seen Mr Appell at the beginning of the months. We talked about you and he forwarded your letter to me.

The theory of continued fractions opens door to a large field of research but where it is not always that easy to figure out in advance the difficulties that will be encountered. It is thus with a pinch of salt that I will recommend  you to read a detailed study of the generalization of continued fractions that should interest you. I only touched upon this topic in a paper that was published few years ago in the journal of Mr. Jordan, which I am happy to enclose to this letter. You will find inside a reference to a paper from Mr. Hermite where he deals with the same topic. In this last paper, Mr. Hermite explains the method of associated polynomials using [...]  integrals: it would be for sure that interesting to go further into the relationship between the calculations using integrals and the recursive rules of the theory of continued fractions.



The result that you said you found seems to be really outstanding but there might be numerous exceptions to it. I would be very happy to read your demonstration when it will be published.
[...] }
\end{quotation}
   \end{minipage}
\vspace{0.8cm}
Moreover, few days later, in a letter\footnote{Collection of letters received by Robert de Montessus de Ballore during the period 1897-1937, from Appell.} sent December the third, 1901, Paul Appell wrote :
\vspace{0.8cm}

\begin{minipage}[c]{.46\linewidth}
\begin{quotation}
{\it Vos r\'esultats me semblent int\'eressants et je suis d'avis que vous les indiquiez [...] dans une note \`a la Soci\'et\'e Math\'ematique [...] }
\end{quotation}
\end{minipage}
\begin{minipage}[c]{.46\linewidth}
\begin{quotation}
{\it 
I think that your results are interesting and you should mention them  [...] in a note to the Soci\'et\'e Math\'ematique [...]
}
\end{quotation}
\end{minipage}
\vspace{0.8cm}

Henri Pad\'e was then professor at the University of Poitiers.


Furthermore, notes from Robert de Montessus in a file named {\it Continued fractions}\footnote{Collection of documents of Robert de Montessus during the period 1897-1937.} enable us to know which work related to continued fraction he consulted at the beginning of his research :
\vspace{0.8cm}

\begin{minipage}[c]{.46\linewidth}
\begin{quotation}
{\it 

Note du 29 Mars 1901

Mr. Hermite gave the expression of the  general convergents of the function $\exp$ [...]}
\end{quotation}
\end{minipage}
\begin{minipage}[c]{.46\linewidth}
\begin{quotation}
\begin{it}
Note of March the 29\up{th}, 1901

Mr. Hermite gave the expression of the  general convergents of the function $\exp$ [...]
\end{it}
\end{quotation}
\end{minipage}
\vspace{0.8cm}
and then a little bit further in this note :
\vspace{0.8cm}

\begin{minipage}[c]{.46\linewidth}
\begin{quotation}
\begin{it}

Sur la fraction de  Stieltjes 19 Avril [...]
\end{it}
\end{quotation}
\end{minipage}
\begin{minipage}[c]{.46\linewidth}
\begin{quotation}
\begin{it}
About the fraction of Stieltjes, April the 19\up{th} [...]
\end{it}
\end{quotation}
\end{minipage}
\vspace{0.8cm}

In the margin of the note are also written the names of  {\it Hurwitz} and {\it Minkowski}. Still in the same file is the article of Salvator Pincherle {\it Sur la g\'en\'eralisation de syst\`emes r\'ecurrents au moyen d'une \'equation diff\'erentielle} published in 1892 in Acta Mathematica. Alfred Pringsheim gave to Robert de Montessus several references for articles in a letter on February the fourth, 1902. 


At the same time, Robert de Montessus was in touch with Eug\`ene Fabry who worked on Taylor series. Fabry sent two letters to Robert de Montessus on October the 12\up{th}, 1901 and January the 7\up{th}, 1902. Emile Borel described and explained the importance of the work of Fabry in his obituary notice read at the Academy of Science on October the 23\up{th}, 1944\footnote{Weekly reports of the meeting of the Academy of Science, 1944, Biblioth\`eque Nationale of France, Gallica.}.  Indeed this notice shows a thorough overview of Fabry's research. In particular, Emile Borel wrote :

\begin{quotation}
{\it Eug\`ene Fabry dedicated then his research to an important and difficult problem that has been already tackled by a lot of mathematicians, such as the our distinguished fellow member Jacques Hadamard. It relates to looking for and staying the singularities located on the convergence circle of a Taylor development that defines an analytical function inside the circle. The brilliant idea of  Fabry was to substitute the study of the complete sequence of the coefficients of the Taylor series by the study of the partial sequences of coefficients extracted from this complete sequence.

 }
\end{quotation}
The two letters from Fabry to  Robert de Montessus are indeed related to Taylor series and their singularities :
\vspace{0.8cm}

\begin{minipage}[c]{.46\linewidth}
\begin{quotation}
(letter of October the 12\up{th}, 1901)
\begin{it}
Le th\'eor\`eme que vous m'indiquez ne peut pas \^etre exact sous la forme la plus g\'en\'erale. Il est possible que le second \'enonc\'e soit exact mais le premier ne l'est pas. Si la s\'erie n'a qu'un point singulier sur la circonf\'erence de convergence, on ne peut affirmer que $\frac{s_{n}}{s_{n+1}}$ ait une limite. Cela r\'esulte des th\'eor\`emes que j'ai indiqu\'e () dans les Acta Math (Tome 22, page 86). J'ai en effet montr\'e qu'il existe des s\'eries incompl\`etes n'ayant qu'un seul point singulier sur la circonf\'erence de convergence. $\frac{s_{n}}{s_{n+1}}$ a alors des valeurs nulles et infinie, et ne peut avoir aucune limite. L'exemple que j'ai donn\'e est le suivant :
$$
\sum x^{n}e^{n\left\lbrack  -1+\cos(L n)^{\theta}\right\rbrack}\quad 0<\theta<1
$$
[...]

En r\'esum\'e, je ne peux pas vous donner de r\'eponse absolument pr\'ecise sur l'exactitude du th\'eor\`eme que vous \'enoncez ; mais je ne serais pas \'etonn\'e qu'il soit exact en prenant le second \'enonc\'e.

Cette question me parait tr\`es int\'eressante et doit conduire \`a des r\'esultats importants.
\end{it}
\end{quotation}
\end{minipage}
\begin{minipage}[c]{.46\linewidth}
\begin{quotation}
(letter of October the 12\up{th}, 1901)
\begin{it}
The theorem you mentioned cannot be correct if it is generalized. It is possible that the second statement is, but not the first one. If the sequence has only one singularity on the circle of convergence, then one cannot affirm that $\frac{s_{n}}{s_{n+1}}$ has a limit. This arises from the theorems I cited [...] in Acta Math (volume 22, page 86). I indeed proved that one can find an incomplete sequence with only one singularity on the circle of convergence. $\frac{s_{n}}{s_{n+1}}$ has then a null or infinite value and cannot have a limit. The example I gave is the following:

$$
\sum x^{n}e^{n\left\lbrack  -1+\cos(L n)^{\theta}\right\rbrack}\quad 0<\theta<1
$$
[...]

To sum up, I cannot give you a very precise answer on the exactness of the theorem you state; however I would not be surprise if it would be true with the second statement.


This issue is very interesting and should lead to important results.

\end{it}
\end{quotation}
\end{minipage}
\vspace{0.8cm}

and
\vspace{0.8cm}

\begin{minipage}[c]{.46\linewidth}
\begin{quotation}
(letter of January the 7\up{th}, 1902)
\begin{it}
[...] Si vous voulez avoir une id\'ee des travaux publi\'es sur la s\'erie de Taylor, vous pourrez trouver des renseignements tr\`es complets dans un petit trait\'e publi\'e par M. Hadamard au mois de mai dernier sur `` la s\'erie de Taylor et son prolongement analytique'' dans la collection Scienta (...).
\end{it}
\end{quotation}
\end{minipage}
\begin{minipage}[c]{.46\linewidth}
\begin{quotation}
(letter of January the 7\up{th}, 1902)
\begin{it}
[...]
If you want to have an idea of the published works on Taylor series,  you can find very complete pieces of information in the little treaty published by Mr. Hadamard last month on "la s\'erie de Taylor et son prolongement analytique'' in the edition Scienta [...].

\end{it}
\end{quotation}
\end{minipage}

This second letter by Fabry enlights the first one. Fabry gave to Robert de Montessus the book of Hadamard as a new reference. Indeed Hadamard took interest in \cite{Hadam2} in the singularities of a function that can be developed in a power series around zero, with a development 
$$
a_{0}+a_{1}x+\cdots +a_{m}x^{m}+\cdots
$$
 having a radius $R$ and a circle of convergence $C$. As did Lecornu\footnote{ Lecornu stated in 1887 in a note to the Comptes Rendus de l'Acad\'emie des Sciences that the existence of a limit for the ratio  $\displaystyle{\frac{a_{m}}{a_{m+1}}}$ would imply that this limit is the unique singularity of the function on  $C$.}, Hadamard considered the ratio  $\displaystyle{\frac{a_{m}}{a_{m+1}}}$ and precised on the page 19 of \cite{Hadam2} :

\begin{quotation}
\begin{it}
Indeed, on the one hand, as we already mentioned, $\displaystyle{\frac{a_{m}}{a_{m+1}}}$ has in general no limit. On the other hand, there could be several singularities on $C$. Then one can see three different meaning in the statement:

\begin{enumerate}
\item If $x_{0}$ is the only singularity on  $C$, the ratio $\displaystyle{\frac{a_{m}}{a_{m+1}}}$ has $x_{0}$ for limit.
\item If  $\displaystyle{\frac{a_{m}}{a_{m+1}}}$ has $x_{0}$ for limit, then the point with affix $x_{0}$ is the only singularity of the function on $C$.

\item If $\displaystyle{\frac{a_{m}}{a_{m+1}}}$ has $x_{0}$ for limit, then the point with affix $x_{0}$ is a singularity of the function.

\end{enumerate}
\end{it}
\end{quotation} 
Jacques Hadamard showed indeed that the first two statements were wrong, contrary to the third one which has been proven by Fabry.


Robert de Montessus moved closer to Pad\'e. Indeed, two other letters by Henri Pad\'e\footnote{Collection of letters received by Robert de Montessus de Ballore during the period1897-1937, from Pad\'e.}  show that Robert de Montessus shared his results with him. Henri Pad\'e\footnote{We should keep in mind that Paul Appell was among the jury of Henri Pad\'e's doctoral defense. He was also one of his professors in Paris.} was at that time professor at the University of Poitiers \cite{Brez1} but we think that Robert de Montessus and him could have met during the international congress of mathematicians that was held in Paris in August 1900. Robert de Montessus got a signed copy of the announcement made by  Henri Pad\'e during this congress, entitled {\it Aper\c cu sur les d\'eveloppements r\'ecents de la th\'eorie des fractions continues.}. Here we should highlight a point related to the vocabulary used : Pad\'e talked of {\it fraction approch\'ee} in his articles to refer to the approximants that have his name now. Van Vleck used the English word  {\it approximant} in 1903 in one of his articles, cited below : it might be this change to the English word that made the French word to be forgotten.


The following excerpts of these two letters underline the mathematical communication between the two men:

\vspace{0.8cm}

\begin{minipage}[c]{.46\linewidth}
\begin{quotation}
(lettre du 17/12/1901)
\begin{it}

Je ne vois que des compliments \`a vous adresser, et la communication que vous proposez de faire \`a la soci\'et\'e math\'ematiques me semble devoir \^etre des plus int\'eressantes.

Vous consid\'erez, d'apr\`es votre lettre, les r\'eduites qui forment la \og ligne horizontale de rang $p$\fg'. Je pense qu'il s'agit des r\'eduites pour lesquelles le degr\'e du \underline {d\'enominateur} est $p$. [...]

Enfin, vous savez que dans mes derniers travaux j'ai chang\'e cette notation et la notion de tableau [...]

Dans cette situation, les fractions qui convergent vers la fonction seraient toutes celles qui correspondent \`a une m\^eme valeur de $\mu$ [...]

\end{it}
\end{quotation}
\end{minipage}
\begin{minipage}[c]{.46\linewidth}
\begin{quotation}
(letter of December the 17\up{th}, 1901)
\begin{it}

I have nothing but congratulations to address to you, and the announcement that you offer to do at the mathematical society seems to be very interesting.


You consider, according to your letter, convergents that form the "horizontal line of rank  $p$''. I think there are actually convergents for which the degree of the \underline {denominator} is $p$. [...]


Finally, you know that in my last work, I changed this notation and the notion of table [...]


In this situation, the fractions that converge toward the function should be those that have the same value of $\mu$ [...]

\end{it}
\end{quotation}
\end{minipage}
\vspace{0.8cm}

Henri Pad\'e suggested to Robert de Montessus to use new notations but the last one will not follow his advice. Robert de Montessus wrote the convergents $(n,p)$ by $\displaystyle{\frac{U_{p}^{n}}{V_{n}^{p}}}$ with $U_{p}^{n}$ a polynomial of degree $n$  and $V_{n}^{p}$ a polynomial of degree $p$ (in the case of a normal table). This convergents corresponds to the Pad\'e approximant $\left\lbrack n,p\right\rbrack$ of the function $f(x)$.

\vspace{0.8cm}

\begin{minipage}[c]{.46\linewidth}
\begin{quotation}
(lettre du 20/6/1902)
\begin{it}
Pour ce qui touche le num\'ero 1 de votre note, sur lequel vous me demandez particuli\`erement mon avis, je vous fais confidence que c'est un point que j'ai d\'ej\`a cherch\'e \`a \'elucider moi-m\^eme. J'ai pr\'esent\'e \`a ce sujet une note \`a l'acad\'emie [...]

Cette co\"incidence n'a pas lieu de vous \'etonner, n'est-ce-pas, et nous ne pouvons pas esp\'erer, que travaillant les m\^emes questions, il n'arrivera jamais que nous parvenions aux m\^emes r\'esultats.
\end{it}
\end{quotation}
\end{minipage}
\begin{minipage}[c]{.46\linewidth}
\begin{quotation}
(Letter of June the 20\up{th}, 1902)
\begin{it}
For everything related to the number 1 of your note, on which you asked for my opinion, I have to confess to you that this a point that I already tried to understand myself. I wrote a note on this topic to the Academy [...]


This coincidence should not surprise you and we could not hope that working on the same issues, we would never end up to the same results.

\end{it}
\end{quotation}
\end{minipage}
\vspace{0.8cm}

Henri Pad\'e referred to his article published in 1902, {\it Recherches nouvelles sur la distribution des fractions approch\'ees d'une fonction}, Annales Scientifiques de l'ENS, volume 19, 1902, pages 153-189. The last sentence of  the excerpt should be put in parallel with the controversy of the priority between Robert de Montessus and Henri Pad\'e that will be further developed later.



\subsection{The theorem}


What is the content of the article published in 1902 ? On which mathematical elements did Robert de Montessus build his demonstration ?

Robert de Montessus started to publish in 1902 articles where he tackled the issue of the convergence of algebraic continued fractions. Thus, before the publication of the theorem in the Bulletin de la Soci\'et\'e Math\'ematique de France, a note from Robert de Montessus was presented at the Academy of Science by Paul Appell on June the 23\up{rd}, 1902. Robert de Montessus got inspiration from a study by Laguerre \cite{Lag3} on the development in continued fractions of the function $\displaystyle{\left(\frac{x+1}{x-1}\right)^{\omega}}$ \cite{RMB1}. He proved the convergence of the continued fraction obtained by Laguerre outside the interval $\left\lbrack -1,1\right\rbrack$.


The theorem of Robert de Montessus  was thus published in the  Bulletin de la Soci\'et\'e Math\'ematique de France in 1902. The title of the article was not precise : {\it Sur les fractions continues alg\'ebriques}. The article was short, with 9 pages. The writing was quite dynamic since the author went directly to the main points. The way how Robert de Montessus has written his article might be one of the reasons why his results has had and still had a special echo. Robert de Montessus derived his theorem as well in his doctoral dissertation : there the writing is even denser. This can be explained by the paragraphs before the result where Robert de Montessus developed with lots of details the links between continued fractions and Pad\'e approximants.


How did  Robert de Montessus build his article ? He started from the power series $y$
$$
y=s_{0}+s_{1}x+\cdots +s_{h}x^{h}+\cdots\quad  (s_{0}\neq 0), 
$$
then reminded the reader several  {\it notions from Pad\'e}: Pad\'e approximants table ; normality of the table and  the necessary condition to have the normality ; link between Pad\'e approximants, seen as convergents, and continued fraction. Robert de Montessus explained that the study of a sequence of fractions well chosen in  Pad\'e table, is similar to the study of a continued fraction. The study of the convergence of the continued fraction, which is the convergence of the sequence of its convergents, is actually the study of the convergence of the series associated to this sequence.


Robert de Montessus imposed some conditions on the degrees of the numerators and denominators of the fractions from Pad\'e table to simplify the problem : if every fraction is {\it in advance of the previous one}, that is to say that $\displaystyle{\frac{U_{p}^{n}}{V_{n}^{p}}}$ is in advance of  $\displaystyle{\frac{U_{q}^{m}}{V_{m}^{q}}}$ if $p+n>m+q$, then the convergents of the sequence have to be all consecutive. This means that if 
$$
\frac{U_{i}}{V_{i}}=\frac{U_{p}^{n}}{V_{n}^{p}},\  \frac{U_{i+1}}{V_{i+1}}=\frac{U_{q}^{m}}{V_{m}^{q}}
$$ 
then $p+n+1=q+m$ or $p+n+2=q+m$.

Robert de Montessus considered only the particular case of the rows in the table, namely the case of a sequence:

$$
\displaystyle{\frac{U_{p}^{0}}{V_{0}^{p}}},\ \displaystyle{\frac{U_{p}^{1}}{V_{1}^{p}}},\ \ldots,\displaystyle{\frac{U_{p}^{n}}{V_{n}^{p}}},\ \ldots.
$$

Then he tackled the two major issues raised in our introduction which are the rebuilding of the function given by the series and its possible analytical continuation. Robert de Montessus considered a meromorphic function, analytical at the origin with a development given by $y$. His proof used the results found by Hadamard \cite{Hadam1} on the asymptotical behavior of polynomials called {\it polynomials of Hadamard}\footnote{P. Henrici in \cite{Henrici1}, pages 622-633, uses these terms.} and of determinants of Hankel $H_{p}^{m}$ associated to the series $y$. The general formulation of the determinants of Hankel is:

$$
H_{p}^{m}=\left\vert
\begin{array}{lrrl}
s_{m}&s_{m+1}&\cdots&s_{m+p-1}\\
s_{m+1}&s_{m+2}&\cdots&s_{m+p}\\
\cdots&\cdots&\cdots&\cdots \\
s_{m+p-1}&s_{m+p}&\cdots&s_{m+2p-2}
\end{array}\right\vert.
$$

They naturally appear when one transform the coefficients of the denominator of the approximate (or of the convergents) in a linear system, with a condition like:

$$
f(z)-\left\lbrack L/M\right \rbrack(z)=O\left(z^{L+M+1}\right)
$$

As for the polynomials of Hadamard, they have been invented by  Carl Jabobi\footnote{Jacobi, C. G. J.,\"Uber die Darstellung einer Reihe gegebner Werthe durch eine gebrochne rationale Function, J. Reine Angew. Math.  30, 127-156.} in 1846. Their general formulation is:

$$
\frac{\left\vert
\begin{array}{lrrll}
s_{m}&s_{m+1}&\cdots&s_{m+p-1}&1\\
s_{m+1}&s_{m+2}&\cdots&s_{m+p}&u\\
\cdots&\cdots&\cdots&\cdots&\vdots \\
s_{m+p-1}&s_{m+p}&\cdots&s_{m+2p-2}&u^{p-1}\\
s_{m+p}&s_{m+p}&\cdots&s_{m+2p-1}&u^{p}
\end{array}\right\vert}
{\left\vert
\begin{array}{lrrl}
s_{m}&s_{m+1}&\cdots&s_{m+p-1}\\
s_{m+1}&s_{m+2}&\cdots&s_{m+p}\\
\cdots&\cdots&\cdots&\cdots \\
s_{m+p-1}&s_{m+p}&\cdots&s_{m+2p-2}
\end{array}\right\vert}.
$$

Indeed, the denominators of the convergents, $V_{n}^{p}$, are nothing else but the Hadamard polynomials of the same degree $p$. Thus, Robert de Montessus raised the hypothesis on the poles 
$$
\vert \alpha_{1}\vert\leq \vert \alpha_{2}\vert\leq \cdots \leq \vert \alpha_{p}\vert<\vert \alpha_{p+1}\vert\leq \vert \alpha_{p+2}\vert\leq\cdots
$$ 
and deduced that the $V_{n}^{p}$ approach a polynomial of degree $p$ whose roots are the poles. Here stands the key point of the result:  Robert de Montessus connected the denominators of the particular Pad\'e approximants that he considered to the result found by Jacques Hadamard page 41 of \cite{Hadam2}. Robert de Montessus continued his reasoning on the series related to the continued fraction:

$$
\frac{U_{p}^{0}}{V_{0}^{p}}+\left( \frac{U_{p}^{1}}{V_{1}^{p}}-\frac{U_{p}^{0}}{V_{0}^{p}}\right)+
\left( \frac{U_{p}^{2}}{V_{2}^{p}}-\frac{U_{p}^{1}}{V_{1}^{p}}\right)+\cdots
$$
It is a series of rational fractions with a known asymptotical behavior of their denominators. This enables us to study the convergence of a series expressed without these denominators. Once it is proved convergent, one should prove that the limit is a meromorphic function. Robert de Montessus thus derivated the following theorem\footnote{The statement of the theorem is given in the conclusion of this article.} :

\begin{quotation}
\begin{it}
[...] given a Taylor series representing a function $f(x)$ whose $p$ poles the closest to the origin are inside the circle $(C)$, this circle being itself inside the next poles, each multiple pole counted as many as simple poles as  it exists unities in his degree, the continued fraction deduced from the row of rank $p$ of Pad\'e table, this table constituted of normal convergents represents a function $f(x)$ in a circle of radius $\vert \alpha_{p+1}\vert$, where $\alpha_{p+1}$ is the affix of the pole the closest to the origin among those outside the circle  $(C)$.

\end{it}
\end{quotation}

Seventy years later, E.B. Saff \cite{Saff_1} stated the theorem of Robert de Montessus as:
\begin{quotation}
\begin{it}
Let $f(z)$ be analytic at $z=0$ and meromorphic with precisely $\nu$ poles (multiplicity counted) in the disk $\vert z\vert <\tau$. Let $D$ the domain obtained from $\vert z\vert <\tau$ by deleting the $\nu$ poles of $f(z)$. Then, for all $n$ sufficiently large, there exits a unique rational function $R_{\mu,\nu}$, of type $(\mu,\nu)$, which interpolates to $f(z)$ in the point $z=0$ considered of multiplicity $\eta+\nu+1$. Each $R_{\eta,\nu}$ has precisely $\nu$ finite poles and, as $n\to\infty$, these pole approach, respectively, the $\nu$ poles of $f(z)$ in $\vert z\vert <\tau$. The sequence $R_{\mu,\nu}$ converges throughout $D$ to $f(z)$, uniformly on any compact subset of $D$.
\end{it}
\end{quotation}
He did not mention continued fractions in the article : interpolation and approximation notions were highlighted instead. Saff referred to the article of J.L. Walsh\footnote{J. L. Walsh, The Convergence of Sequences of Rational Functions of Best Approximation. II, Transactions of the American Mathematical Society, Vol. 116 (Apr., 1965), pp. 227-237.} in a previous article\footnote{Saff E.B., On the row convergence of the Walsh array of meromorphic functions, Transactions of the American Mathematical Society, vil. 146, 1969.}. In the introduction of this article, J.L. Walsh wrote :

\begin{quotation}
\begin{it}
The $W_{n\mu}$ form a table of double entry (...) known as the Walsh array which is similar in form and properties to the table of Pad\'e. Indeed, J.L. Walsh has for the rows of this array  established (...) the following analogue of the important result (...) of Montessus de Ballore concerning the convergence of the rows of the Pad\'e table.
\end{it}
\end{quotation}

\subsection{The work of Robert de Montessus and of Henri Pad\'e}


As we mentioned it before, Henri Pad\'e and Robert de Montessus were in touch since the beginning of the 20\up{th} century. This did not lead to any collaboration excepted the request of off-print articles\footnote{Collection Robert de Montessus, off-print articles of Pad\'e, with the author's signature on it.}. Let us give some quantitative elements on the work of these two authors in the field of continued fractions. Henri Pad\'e published 28 articles between 1890 and 1907\footnote{Claude Brezinski in \cite{Brez1} counted 42 articles in mathematics.} whereas Robert de Montessus published  9 articles\footnote{ Robert de Montessus published a total of 57 articles but also 13 books and several papers in journals such as  {\it L'interm\'ediaire des math\'ematiciens} (founded by C.A. Laisant), {\it L'enseignement math\'ematique} (founded by C.A Laisant et H. Fehr), {\it La revue du mois} (founded by E. Borel and his wife the writer Camille Marbo, daughter of P. Appell).} on algebraic continued fractions between 1902 and 1909. Robert de Montessus has preferentially worked on the convergence whereas Pad\'e more on the creation of his approximants and on the methods to obtain them. Since 1908,  Pad\'e occupied several jobs of state superintendent of education ; he thus stopped his research. However we do not know why Robert de Montessus did not published any other work on the algebraic continued fractions after the year 1909.


Robert de Montessus referred several times as expected to the work of Pad\'e. However, how did Henri Pad\'e received the results of Robert de Montessus ? What could the reasons for Robert de Montessus to send a letter\footnote{Letter of Emile Picard written on the first of January 1908, collection Robert de Montessus.} to  Emile Picard, director of the Annales Scientifiques of the Ecole Normale Sup\'erieure to protest at the end of the year 1907 ? The starting point of the controversy is the last article of Pad\'e on algebraic continued fractions: {\it Recherches sur la convergence des d\'eveloppements en fractions continues d'une certaine cat\'egorie de fonctions} \cite{Pade1}. Even if Robert de Montessus talked of {\it matter of priority between him and me} in the excerpt of the letter published in the Annales \cite{RMB2}, the reaction of Robert de Montessus is related to the convergence of the development in continued fractions of the function

 $$
\displaystyle{\left(\frac{x-1}{x+1}\right)^{\omega}}.
$$
Finally, was there a really reason for this controversy ? Henri Pad\'e was objectively right when he wrote in his article of 1907, page 358 :
\begin{quotation}
\begin{it}
To conclude, with the convergence of the series equivalent to the continued fraction supposed to be valid, how is it that its sum, i.e. the value of the continued fraction, is equal to the function $\displaystyle{\left(\frac{x-1}{x+1}\right)^{\omega}}$? \\
No indication of the derivation of this major point is furnished in the notes of Mr. de Montessus.

\end{it}
\end{quotation}
But Henri Pad\'e leant on three articles of Robert de Montessus published in 1902, 1903 and 1904, forgetting his dissertation of 1905 and a note to the Comptes rendus de l'Acad\'emie des Sciences of Mai the 29\up{th}, 1905, where Robert de Montessus proved that the limit of the convergents is indeed the function. The function above is a particular case of functions resolving differential equations such as

$$
(az+b)(cz+d)Z'(z)=(pz+q)Z(z)+P(z)
$$
where $a$, $b$, $c$, $d$, $p$ and $q$ are constants and $P$ a polynomial. Robert de Montessus resumed the work of Laguerre on the development in continued fraction of such a function. Then Pad\'e again considered this issue but using any kind of continued fraction from Pad\'e table, and he was looking for a integral expression between the function and the convergent, as he already did with the function exponential. This view is always present in the work of Pad\'e, in the contrary to Montessus. The more obvious example is his article of 1902 where there was no {\it formula} giving the rest. Finally, is not the reaction of Robert de Montessus a little bit tough ? In the end, the works of Robert de Montessus and of Pad\'e are complementary to each other.


\section{Circulation of the theorem}


\subsection{Before 1914}

Three factors promoted the rapid circulation of the theorem before 1914. The first one is the references to the results of Robert de Montessus in the work of Pad\'e, Van Vleck and N\"orlund. It is worth reminding that Robert de Montessus shared his result with Henri Pad\'e at the beginning of the year 1902. Van Vleck mentioned the result in a seminar in Boston in 1903. Indeed, he wrote in {\it Pad\'e's Table of Approximants and its Applications.} \cite{VanV} :

\begin{quotation}
\begin{it}
In investigating the convergence of the horizontal lines the first
case to be considered is naturally that of a function having a number
of poles and no other singularities within a prescribed distance of
the origin. It is just this case that Montessus [33, a] has exam
ined very recently. Some of you may recall that four years ago in
the Cambridge colloquium Professor Osgood took Hadamard's
thesis as the basis of one of his lectures. This notable thesis is
devoted chiefly to series defining functions with polar singularities.
Montessus builds upon this thesis and applies it to a table possess
ing a normal character. Although his proof is subject to this
limitation, his conclusion is nevertheless valid when the table is
not normal, as I shall show in some subsequent paper.
\end{it}
\end{quotation}

Van Vleck clarified what we have previously pointed out, that is to say that Robert de Montessus used the results of Hadamard on the polar singularities. On top of that, Van Vleck mentioned a new aspect, namely that the theorem is still accurate even if Pad\'e table is not normal. Moreover, N.E. N\"orlund sent the following\footnote{Collection of letters received by Robert de Montessus de Ballore, period 1897-1937, letter from N.E. N\"orlund sent on March the 29\up{th}, 1910. } to Robert de Montessus in 1910 :

\vspace{0.8cm}

\begin{minipage}[c]{.46\linewidth}
\begin{quotation}
\begin{it}
Les \og Rendiconti di Palermo\fg ne se trouvent \`a aucune biblioth\`eque publique de Copenhague, mais j'ai obtenu aujourd'hui vos th\`eses.

Je suis heureux maintenant de pouvoir citer votre m\'emoire en reconnaissant votre priorit\'e. Le m\'emoire dont j'ai eu l'honneur de vous envoyer un tirage \`a part, ne para\^itra que dans le tome 34 des Acta Mathematica en \underline{1911} [...]

\end{it}
\end{quotation}
\end{minipage}
\begin{minipage}[c]{.46\linewidth}
\begin{quotation}
\begin{it}
I cannot find the \og Rendiconti di Palermo\fg in any public library  of Copenhagen, but I managed to get your dissertations today.

I am glad now to be able to cite your thesis and to acknowledge your priority. The thesis I have the pleasure to send to you an off-print will be published in the volume 34 of Acta Mathematica in \underline{1911} [...]


\end{it}
\end{quotation}
\end{minipage}
\vspace{0.8cm}

This dissertation of more than a hundred pages is entitled {\it Fractions continues et diff\'erences r\'eciproques }.


The second factor of the circulation is given by Robert de Montessus himself. Indeed, after his defense in 1905, he got in touch with Mittag-Leffler\footnote{Collection of letters received by Robert de Montessus de Ballore, period 1897-1937, lettre sent by Mittag-Leffler on December the 30\up{th}, 1905. Robert de Montessus wrote on the letter: {\it dissertation not sent}.}. His results on continued fractions were published in Acta Mathematica, the journal owned by Mittag-Leffler in 1909. His dissertation was also published in a foreign magazine in 1905, Rendiconti di Palermo.

\begin{table}[h]

\begin{center}
\begin{tabular}{|c|c|p{3.5cm}|p{3.5cm}|}
\hline
year&author&title&publisher\\
\hline
1905&Robert de Montessus&Sur les fractions continues alg\'ebriques& volume XIX, pages 185-257,  Rendiconti del Circolo Matematico di Palermo.\\
\hline
1909&Montessus de Ballore R.&Les fractions continues alg\'ebriques&Acta mathematica,volume 32, pages 257-282\\
\hline

\end{tabular}
\end{center}
\caption{publications by Robert de Montessus in 1905 and 1909}

\end{table}

Finally, the academic network contributed in the usability of the result. The award of a part of the Grand Prix de l'Acad\'emie des Sciences in 1906\footnote{The Grand Prix of 1906 corresponded in fact in a throw in of a Grand Prix proposed in 1904.} has also certainly boosted the larger circulation of his result. 


The table underneath gives the articles that contain a reference to the theorem of Robert de Montessus. One can find the article of N\"orlund that we already mentioned. As for the work of  Oskar Perron, it is the first version of the paper that have been published several times after the first world war.

\begin{table}[h]

\begin{center}
\begin{tabular}{|c|c|p{3.5cm}|p{3.5cm}|}
\hline
year&author&title&publisher\\
\hline
1903&Van Vleck E.B.&Selected topics in the theory of divergent series and of continued fractions & Lectures on mathematics, AMS\\
\hline
1907&Pad\'e H.&Sur la g\'en\'eralisation des formules de Sylvester relatives aux fonctions qui se pr\'esentent dans l'application du th\'eor\`eme de Sturm, et sur la convergence des r\'eduites d'une fraction rationnelle&Annales scientifiques de l'ENS, third series, volume 24, pages 519-534\\
\hline
1910& Watson G. N.&The solution of a certain transcendental equation& [J] Lond. M. S. Proc. (2) 8, 162-177\\
\hline
1911&N\"orlund, N. E.&Fractions continues et diff\'erences r\'eciproques&Acta Math. 34, 1-108\\
\hline
1913&Perron O.&Die Lehre von den Kettenbruchen&Leipzig und Berlin, Druck und Verlag von B.G. Teubner\\
\hline
\end{tabular}
\end{center}
\caption{Authors mentioning the theorem in the years 1903-1913}

\end{table}

\subsection{In between the two world wars}
The result was not forgotten after the first world war. Several English or American authors took interest in the theorem. As an example, the English mathematician R. Wilson contacted\footnote{Collection of letters received by Robert de Montessus de Ballore, period 1897-1937, letters from Wilson sent on October the 19\up{th} and 30\up{th}, 1923. } Robert de Montessus in 1923:

\begin{quotation}
(letter of October the 19\up{th}, 1923)
\begin{it}
I am interested in the developement of M. Pad\'es work on the representation of a function by means of a continued fractions, and find that you have contributed much recent to this subject [...]
\end{it}
\end{quotation}
and
\begin{quotation}
(letter of October the 30\up{th}, 1923)
\begin{it}
Your important work on the continued fractions, and the singularities of this function, is special interest to me [...] 
\end{it}
\end{quotation}

The American mathematician J.L. Walsh referred to the theorem in a review of the A.M.S. published in 1935. It appears that he was a contributor to the circulation of the result. Indeed, he was in Paris in 1920-1921, where he worked with Paul Montel. He defended his thesis in 1927 under the supervision of  Maxime B\^ocher. We shall mention here that Maxime B\^ocher stayed in Paris in the years 1913-1914 at the occasion of an exchange between the Universities of Harvard and Paris. Maxime B\^ocher had F. Klein for his thesis supervision. Yet, E.B. Van Vleck who we mentioned several times already, defended as well his thesis\footnote{In 1893.} under the supervision of F. Klein. Then he supervised himself H.S. Wall in 1927, who mentioned then the theorem of Robert de Montessus de Ballore. In the sixties\footnote{To be more precise, in 1964, 1965 and 1967.} J.L. Walsh referred to the theorem of Robert de Montessus in a series of articles. Thus, J.L. Walsh induced the circulation of the theorem \footnote{Maybe he is not the only one?} after the second world war.


The following table gives the articles containing a reference to the theorem of Robert de Montessus and published during the period 1924-1940. N\"orlund quoted the theorem again in 1924; Perron as well in the second edition of his book on continued fractions.


\begin{table}[h]

\begin{center}
\begin{tabular}{|c|c|p{3.5cm}|p{3.5cm}|}
\hline
year&author&title&publisher\\
\hline
1924&N\"orlund N. E.&Vorlesungen uber Differenzenrechnung&Berlin Verlag von Julius Springer\\
\hline
1927&Wilson R.&Divergent continued fractions and polar singularities&Proceedings L. M. S. (2) 26, 159-168\\
\hline
1928&Wilson R.&Divergent continued fractions and polar singularities&Proceedings L. M. S. (2) 27, 497-512; (2) 28, 128-144\\
\hline
1929&Perron O.&Die Lehre von den Kettenbruchen&Second edition, revised. Leipzig and Berlin, Teubner\\
\hline
1935&Walsh J.L.&Interpolation and approximation by rational functions in the complex domain&New York, American Mathematical Society (Amer. Math. Soc. Colloquium Publ. Vol. XX)\\
\hline
1939&Scott W.T., Wall H.S.&Continued Fractions&National Mathematics Magazine, Vol. 13, No. 7, pp. 305-322\\
\hline
1940&Mall J.\footnote{Student of O. Perron}&Beitrag zur Theorie der mehrdimensionalen Pad\'eschen Tafel&[J] Math. Z. 46, 337-349\\
\hline
\end{tabular}
\end{center}

\caption{Author referring to the theorem in the years 1924-1940}


\end{table}


\section{Conclusion}

To conclude this study, it appeared that the result of Robert de Montessus is part of a group of research on the issue of the convergence of algebraic continued fractions led by several mathematicians since the middle of the 19\up{th} century. We have also linked this issue to the problem of the analytical continuation and to other questions in the analysis domain, such as the summation of divergent series or the moment problem. Robert de Montessus summarized in his article of 1902 the work of Pad\'e on the link between taylor series, continued fractions and {\it fractions approch\'ees}, and the work of Hadamard on polar singularities. These are the strength of his work and the reason it has been quickly circulated.


The circulation took place in different ways: by publishing in foreign magazines\footnote{H\'el\`ene Gispert points out that French mathematician did not published in foreign magazines, except Henri Poincar\'e \cite{Gispert}.}, thanks to the will of Robert de Montessus to communicate with the mathematical community of his time, and at last by choosing a relevant topic.


After the second world war, the result of Robert de Montessus tackled already more domains than the continued fractions only, namely interpolation and approximation. The American mathematician J. L. Walsh, who also wrote in French and in German, played an important role. No communications with Robert de Montessus are known.


The position of Robert de Montessus at the editorial board of the Journal de Math\'ematiques Pures et Appliqu\'ees\footnote{Robert de Montessus joined the editorial board in 1917. The journal was then managed by Camille Jordan.}  might have played a part in the circulation of his result as well.


\begin{thebibliography}{1}
\bibitem{site1} http://www.numdam.org/
\bibitem{site2} http://www.emis.de/MATH/JFM/JFM.html
\bibitem{site3} http://www.jstor.org/
\bibitem{site4} http://gallica.bnf.fr/
\bibitem{site5} http://www.ams.org/mathscinet/
\bibitem{site6} http://www.archive.org/
\bibitem{site7} http://gdz.sub.uni-goettingen.de/en/
\bibitem{doc1} Registre AJ/16/5514, Proc\`es verbaux d'examens, Doctorats d'\'etat, du 1/11/1895 au 11/11/1907, num\'ero d'ordre 1187, Archives Nationales, CRAN, Paris.
\bibitem{doc8} Registre AJ/16/5538, rapports de th\`eses Facult\'e des Sciences de Paris, Archives Nationales, CRAN, Paris.
\bibitem{Appell_1} Appell P., Sur les fractions continues p\'eriodiques, Archiv der Mathematik und Physik, 1877.
\bibitem{Appell_2} Appell P., Notice sur les travaux scientifiques de M. Paul Appell, Gauthier-Villars, 1892.
\bibitem{Baker} Baker G.A., Graves-Morris P., Pad\'e Approximants, Cambridge University Press, 1996.
\bibitem{Borel} Borel E., Le\c cons sur s\'eries divergentes, Gauthier-Villars, Paris, 1901.
\bibitem{Borel2} Borel, E., Le\c cons sur les s\'eries divergentes. 2. \'ed. revue et enti\`erement remani\'ee avec le concours de G. Bouligand, Paris: Gauthier-Villars. (Collection de monographies sur la th\'eorie des fonctions). IV, 260 p. (1928).
\bibitem{Brez1} Brezinski Claude, H. Pad\'e, Oeuvres, Librairie scientifique et technique Albert Blanchard, Paris, 1984.
\bibitem{Brez2} Brezinski Claude, History of continued fractions and Pad\'e approximants, Springer Verlag, Berlin, 1991.
\bibitem{Cuyt} Cuyt A., Lubinsky D.S., {\it A de Montessus theorem for multivariate homogeneous Pad\'e approximants.}, Annals of numerical mathematics 4 (1997), pp 217-228.
\bibitem{Dieudonnee} J Dieudonn\'e, Fractions continu\'ees et polyn\^omes orthogonaux dans l'oeuvre de E. N. Laguerre in  : Polynomes orthogonaux et applications, C. Brezinski, A. Draux, A.-P. Magnus, P. Maroni, Lectures Notes in Mathematics, 1171,pp 1-15,  Springer-Verlag, 1985 
\bibitem{Gispert} Gispert H., La France math\'ematique, La soci\'et\'e math\'ematique de France (1870-1914), Cahiers d'histoire et de philosophie des sciences, num\'ero 34, 1991.
\bibitem{Gold_1} Goldstein C., Un th\'eor\`eme de Fermat et ses lecteurs, collection Histoire des Sciences, Presses Universitaire de Vincennes, 1995.
\bibitem{Gut} Gutknecht M.H., Parlett B.N., From $QD$ to $LR$ and $QR$, or, how were the $QD$ and $LR$ algorithms discovered ? ,  IMA J Numer Anal (2011) 31(3): 741-754.
\bibitem{Hadam1} Hadamard J., Essai sur l'\'etude des fonctions donn\'ees par leur d\'eveloppement de Taylor, JMPA, 8, pp 101-186, 1892.
\bibitem{Hadam2} Hadamard J., La s\'erie de Taylor et son prolongement analytique, Scienta Phys. Math., 12, 1901.
\bibitem{Henrici1} Henrici P., Applied and computational complex analysis, volume I, John Wiley and sons, New York, 1988.
\bibitem{Henrici2} Henrici P., Applied and computational complex analysis, volume II, John Wiley and sons, New York, 1991.
\bibitem{Huygens} Huygens C., Oeuvres compl\`etes. Tome XXI. Cosmologie (ed. J.A. Vollgraff). Martinus Nijhoff, Den Haag 1944. 
\bibitem{Jones} Jones W.B., Thron W.J., Continued Fractions, Analytic theory and applications, Encyclopedia of mathematics ans its applications, Addison-Wesley, 1980.
\bibitem{Koch} Von Koch H., Sur un th\'eor\`eme de Stieltjes et sur les fonctions d\'efinies par des fractions continues, Bulletin de la S.M.F., tome 23, pp 33-40, 1895.
\bibitem{Legendre} Legendre A. M., El\'ements de G\'eom\'etrie, avec notes, suivi d'un trait\'e de Trigonom\'etrie, Langlet et Compagnie, Libraires, Bruxelles, 1837.
\bibitem{Lef} Notice sur la vie et l'oeuvre de Robert de Montessus de Ballore, in Dictionnaire biographoque de la Facult\'e des Sciences de Paris 1808-1940, INRP, \`a para\^itre en 2011.
\bibitem{Lag3} Laguerre E., Sur la fonction $\displaystyle{\left(\frac{x+1}{x-1}\right)^{\omega}}$, Bulletin de la SMF, tome 8, pp 36-52, 1880.
\bibitem{Lag1} Laguerre E., Oeuvres compl\`etes, volume 1, sous la direction de C. Hermite, H. Poincar\'e et E. Rouch\'e, Gauthier-Villars, 1898.
\bibitem{Lag2} Laguerre E., Oeuvres compl\`etes, volume 2, sous la direction de C. Hermite, H. Poincar\'e et E. Rouch\'e, Gauthier-Villars, 1905.
\bibitem{Leloup} Leloup J., L'entre-deux-guerres math\'ematique \`a travers les th\`eses soutenues en France, th\`ese de doctorat, Universit\'e Pierre et Marie Curie, 2009.
\bibitem{Maz} Mas'ya V., Shaposhnikova T., Jacques Hadamard, un math\'ematicien universel, EDP Sciences, 2005.
\bibitem{RMB1}De Montessus Robert, Sur les fractions continues alg\'ebriques, Comptes Rendus 134 (1902), 1489-1491. 
\bibitem{RMB4}De Montessus Robert, Sur les fractions continues alg\'ebriques, Bulletin de la SMF, tome 30 (1902), p 28-36.
\bibitem{RMB2}De Montessus Robert, Sur les fractions continues alg\'ebriques (extrait d'une lettre adress\'ee \`a la r\'edaction), Annales Scientifique de l'ENS, 3i\`eme s\'erie, tome 25, 1908, p 195-197
\bibitem{RMB5}De Montessus Robert, Les fractions continues alg\'ebriques, Acta Mathematica, 32, pp 257-281, 1909.
\bibitem{RMB3}De Montessus Robert, notice parue dans le Journal de Math\'ematiques pures et Appliqu\'ees, tome 16, 1937
\bibitem{Padethese} Pad\'e H., Sur la repr\'esentation approch\'ee d'une fonction par des fractions rationnelles, Gauthier-Villars, Paris 1892.
\bibitem{Pade1} Pad\'e H., Recherches sur la convergence des d\'eveloppements en fractions continues d'une certaine cat\'egorie de fonctions, Ann. Sci. Ec. Norm. Sup., 24, pp 341-400, 1907.
\bibitem{Pad\'e_1} Pad\'e H., Sur la repr\'esentation approch\'ee d'une fonction par des fractions rationnelles, Annales de l'Ecole Normale Sup\'erieure, 3 i\`eme s\'erie, tome 9, suppl\'ement, pp 3-93, 1892.
\bibitem{Pad\'e_2} Pad\'e H., Recherches sur la convergence des d\'eveloppements en fractions continues d'une certaine cat\'egorie de fonctions, Ann. Ec. Norm. Sup., 24 (1907), pp 341-400.
\bibitem{Perron} Perron O. {\it Die Lehre von den Kettenbruchen},Teubner, Stuttgart (\'edition de ) 1957.
\bibitem{Pincher} Pincherle S., Notice sur les travaux Acta Mathematica, Vol. 46 (1925), p. 341–362.
\bibitem{Poincare1} Poincar\'e H., Sur un mode nouveau de repr\'esentation g\'eom\'etrique des formes quadratiques d\'efinies ou ind\'efinies, Journal de l'Ecole Polytechnique, tome 28, cahier 47, pp 177-245, 1880.
\bibitem{Poincare2} Poincar\'e H., Sur une g\'en\'eralisation des fractions continues, C.R. Acad. Sci., Paris, 99, (1884), 1014-1016.
\bibitem{Poincare3} Poincar\'e H., Les m\'ethodes nouvelles de la m\'ecanique c\'eleste. Gauthier-Villars, Paris, tome 1,2 et 3, 1892, 1893 et 1899.
\bibitem{Rouse} Histoire des math\'ematiques par W W Rouse Ball ; \'edition fran\c caise revue et augment\'ee, traduite sur la troisi\`eme \'edition anglaise par L Freund, Ball, Walter William Rouse. Auteur; Montessus de Ballore, Robert ; \'editeur scientifique; L Freund, Traduction, Hermann, 1906-1907. -2 vol. (VII-422 p., 271 p.).
\bibitem{Saff_1} E. B. Saff, An extansion of Montessus de Ballore's theorem on the convergence of interpolating rational functions, Journal of Approximation Theory, vol 6, No. 1, July 1972.
\bibitem{Sina}H. Sinaceur, Corps et mod\`eles, Mathesis, Vrin, 1991.
\bibitem{Thom} Thom\'e, L. W., Ueber die Kettenbruchentwicklung des Gaussschen Quotienten $\frac{F(\alpha,\beta+1,\gamma+1,x)}{F(\alpha,\beta,\gamma,x)}$,Reine Angew. Math.  67, 299-309, 1867.
\bibitem{bioVan} Edward Burr Van Vleck, 1863-1943, A Biographical Memoir by Rudolph E. Langer and Mark H. Ingraham, 1957,  National Academy of Sciences
Washington d.c.
\bibitem{VanV} Van Vleck E.B., {\it Selected topics in the theory of divergent series and of continued fractions.}, in Lectures on mathematics
deivered from September 2 to 5, 1903, before members of the american mathematical society
in connexion with the summer meeting hed at the MIT Boston, mass., published for the american mathematical society by the Macmillan company London, 1905. 
\bibitem{Vara} V. S. Varadarajan, Euler and his work on infinite series, Bulletin (New Series) of the American Mathematical Society, Volume 44, Number 4, October 2007, Pages 515-539

\end{thebibliography}
\end{document}